# *A Family of Estimators of Population Mean Using Multiauxiliary Information in Presence of Measurement Errors*

**Jack Allen[1], Housila P. Singh[2], Florentin Smarandache[3]**

[1] School of Accounting and Finance, Griffith University , Gold Coast Campus, Queensland, Australia
[2] School of Statistics, Vikram University, UJJAIN 456010, India
[3] Department of Mathematics, University of New Mexico, Gallup, USA

## Abstract

This paper proposes a family of estimators of population mean using information on several auxiliary variables and analyzes its properties in the presence of measurement errors.

### 1. INTRODUCTION

The discrepancies between the values exactly obtained on the variables under consideration for sampled units and the corresponding true values are termed as measurement errors. In general, standard theory of survey sampling assumes that data collected through surveys are often assumed to be free of measurement or response errors. In reality such a supposition does not hold true and the data may be contaminated with measurement errors due to various reasons; see, e.g., Cochran (1963) and Sukhatme *et al* (1984).

One of the major sources of measurement errors in survey is the nature of variables. This may happen in case of qualitative variables. Simple examples of such variables are intelligence, preference, specific abilities, utility, aggressiveness, tastes, etc. In many sample surveys it is recognized that errors of measurement can also arise from the person being interviewed, from the interviewer, from the supervisor or leader of a team of interviewers, and from the processor who transmits the information from the recorded interview on to the punched cards or tapes that will be analyzed, for instance, see Cochran (1968). Another source of measurement error is when the variable is conceptually well defined but observations can be obtained on some closely related substitutes termed as proxies or surrogates. Such a situation is encountered when one needs to measure the economic status or the level of education of individuals, see Salabh (1997) and Sud and Srivastava (2000). In presence of measurement errors, inferences may be misleading, see Biemer *et al* (1991), Fuller (1995) and Manisha and Singh (2001).

There is today a great deal of research on measurement errors in surveys. An attempt has been made to study the impact of measurement errors on a family of estimators of population mean using multiauxiliary information.

## 2. THE SUGGESTED FAMILY OF ESTIMATORS

Let Y be the study variate and its population mean $\mu_0$ to be estimated using information on p(>1) auxiliary variates $X_1, X_2, ..., X_p$. Further, let the population mean row vector $\underset{\sim}{\mu}' = (\mu_1, \mu_2, \cdots, \mu_p)$ of the vector $\underset{\sim}{X}' = (X_1, X_2, X_p)$. Assume that a simple random sample of size n is drawn from a population, on the study character Y and auxiliary characters $X_1, X_2, ..., X_p$. For the sake of simplicity we assume that the population is infinite. The recorded fallible measurements are given by

$$\begin{aligned} y_j &= Y_j + E_j \\ x_{ij} &= X_{ij} + \eta_{ij}, \quad i = 1,2,\cdots, p; \\ &\qquad\qquad\quad\; j = 1,2,\cdots, n. \end{aligned}$$

where $Y_j$ and $X_{ij}$ are correct values of the characteristics Y and $X_i$ (i=1,2,..., p; j=1,2,..., n).

For the sake of simplicity in exposition, we assume that the error $E_j$'s are stochastic with mean 'zero' and variance $\sigma_{(0)}^2$ and uncorrelated with $Y_j$'s. The errors $\eta_{ij}$ in $x_{ij}$ are distributed independently of each other and of the $X_{ij}$ with mean 'zero' and variance $\sigma_{(i)}^2$ (i=1,2,...,p). Also $E_j$'s and $\eta_{ij}$'s are uncorrelated although $Y_j$'s and $X_{ij}$'s are correlated.

Define

$$\begin{aligned} u_i &= \frac{\bar{x}_i}{\mu_i}, \ (i = 1,2,\cdots, p) \\ u^T &= (u_1, u_2, \cdots u_p)_{1\times p}, \ e^T = (1,1,\cdots,1)_{1\times p} \\ \bar{y} &= \frac{1}{n}\sum_{j=1}^{n} y_j \\ \bar{x}_i &= \frac{1}{n}\sum_{j=1}^{n} x_{ij} \end{aligned}$$

With this background we suggest a family of estimators of $\mu_0$ as

$$\hat{\mu}_g = g(\bar{y}, u^T)$$

(2.1)

where $g\left(\bar{y},u^T\right)$ is a function of $\bar{y},u_1,u_2,\cdots,u_p$ such that

$$g_{\left(u_0,e^T\right)}=\mu_0$$
$$\Rightarrow \left.\frac{\partial g(\cdot)}{\partial \bar{y}}\right|_{\left(\mu_0,e^T\right)}=1$$

and such that it satisfies the following conditions:

1. The function $g\left(\bar{y},u^T\right)$ is continuous and bounded in Q.

2. The first and second order partial derivatives of the function $g\left(\bar{y},u^T\right)$ exist and are continuous and bounded in Q.

To obtain the mean squared error of $\hat{\mu}_g$, we expand the function $g\left(\bar{y},u^T\right)$ about the point ($\mu_0$,$e^T$) in a second order Taylor's series. We get

$$\hat{\mu}_g=g\left(\mu_0,e^T\right)+\left(\bar{y}-\mu_0\right)\left.\frac{\partial g(\cdot)}{\partial \bar{y}}\right|_{\left(\mu_0,e^T\right)}+(u-e)^T g^{(1)}{}_{\left(\mu_0,e^T\right)}$$
$$+\frac{1}{2}\left\{\left(\bar{y}-\mu_0\right)^2\left.\frac{\partial^2 g(\cdot)}{\partial \bar{y}^2}\right|_{\left(\bar{y}^*,u^{*T}\right)}+2\left(\bar{y}-\mu_0\right)(u-e)^T\left.\frac{\partial g^{(1)}(\cdot)}{\partial \bar{y}}\right|_{\left(\bar{y}^*,u^{?T}\right)}\right.$$
$$\left.+(u-e)^T g^{(2)}\left(\bar{y}^*,u^{*T}\right)(u-e)\right\}$$

(2.2)

where

$$\bar{y}^*=\mu_0+\theta\left(\bar{y}-\mu_0\right),\ u^*=e+\theta(u-e),\ 0<\theta<1;\quad g^{(1)}(\cdot)$$

denote the p element column vector of first partial derivatives of g(·) and $g^{(2)}(\cdot)$ denotes a p×p matrix of second partial derivatives of g(·) with respect to u.

Noting that $g(\mu_0,e^T)=\mu_0$, it can be shown that

$$E\left(\hat{\mu}_g\right)=\mu_0+O(n^{-1})$$

(2.3)

which follows that the bias of $\hat{\mu}_g$ is of the order of $n^{-1}$, and hence its contribution to the mean squared error of $\hat{\mu}_g$ will be of the order of $n^{-2}$.

From (2.2), we have to terms of order $n^{-1}$,

$$\begin{aligned}\text{MSE}(\hat{\mu}_g) &= E\left\{(\bar{y}-\mu_0)+(u-e)^T g^{(1)}(\mu_0,e^T)\right\}^2 \\ &= E\left[(\bar{y}=\mu_0)^2 + 2(\bar{y}-\mu_0)(u-e)^T g^{(1)}(\mu_0,e^T)\right. \\ &\qquad \left. + \left(g^{(1)}(\mu_0,e^T)\right)^T (u-e)(u-e)^T \left(g^{(1)}(\mu_0,e^T)\right)\right] \\ &= \frac{1}{n}\left[\mu_0^2\left(C_0^2 + C_{(0)}^2\right) + 2\mu_0 b^T g^{(1)}(\mu_0,e^T) + \left(g^{(1)}(\mu_0,e^T)\right)^T A\left(g^{(1)}(\mu_0,e^T)\right)\right]\end{aligned}$$

(2.4)

where $b^T=(b_1,b_2,\ldots,b_p)$, $b_i,=\rho_{0i}C_0C_i,(i=1,2,\ldots,p)$;

$C_i=\sigma_i/\mu_i$, $C_{(i)}=\sigma_i/\mu_i$, $(i=1,2,\ldots,p)$ and $C_0=\sigma_0/\mu_0$,

$$\text{A} = \begin{bmatrix} C_1^2 + C_{(1)}^2 & \rho_{12}C_1C_2 & \rho_{13}C_1C_3 & \cdots & \rho_{1p}C_1C_p \\ \rho_{12}C_1C_2 & C_2^2 + C_{(2)}^2 & \rho_{23}C_2C_3 & \cdots & \rho_{2p}C_2C_p \\ \rho_{13}C_1C_3 & \rho_{23}C_2C_3 & C_3^2 + C_{(3)}^2 & \cdots & \rho_{3p}C_3C_p \\ \vdots & \vdots & \vdots & & \vdots \\ \rho_{1p}C_1C_p & \rho_{2p}C_2C_p & \rho_{3p}C_3C_p & \cdots & C_p^2 + C_{(p)}^2 \end{bmatrix}_{p\times p}$$

The $\text{MSE}(\hat{\mu}_g)$ at (2.4) is minimized for

$$g^{(1)}(\mu_0,e^T) = -\mu_0 A^{-1} b$$

(2.5)

Thus the resulting minimum MSE of $\hat{\mu}_g$ is given by

$$\text{min.MSE}(\hat{\mu}_g) = \left(\mu_0^2/n\right)\left[C_0^2 + C_{(0)}^2 - b^T A^{-1} b\right]$$

(2.6)

Now we have established the following theorem.

Theorem 2.1 = Up to terms of order $n^{-1}$,

$$MSE(\hat{\mu}_g) \geq \left(\mu_0^2/n\right)\left[C_0^2 + C_{(0)}^2 - b^T A^{-1} b\right]$$

(2.7)

with equality holding if

$$g^{(1)}(\mu_0,e^T) = -\mu_0 A^{-1} b$$

It is to be mentioned that the family of estimators $\hat{\mu}_g$ at (2.1) is very large. The following estimators:

$$\hat{\mu}_g^{(1)} = \bar{y}\sum_{i=1}^{p}\omega_i\left(\frac{\mu_i}{\bar{x}_i}\right); \quad \sum_{i=1}^{p}\omega_i = 1,$$

[Olkin (1958)]

$$\hat{\mu}_g^{(2)} = \bar{y}\sum_{i=1}^{p}\omega_i\left(\frac{\bar{x}_i}{\mu_i}\right), \quad \sum_{i=1}^{p}\omega_i - 1, \text{ [Singh (1967)]}$$

$$\hat{\mu}_g^{(3)} = \bar{y}\frac{\sum_{i=1}^{p}\omega_i\mu_i}{\sum_{i=1}^{p}\omega_i\bar{x}_i}, \quad \sum_{i=1}^{p}\omega_i = 1, \text{ [Shukla (1966) and John (1969)]}$$

$$\hat{\mu}_g^{(4)} = \bar{y}\frac{\sum_{i=1}^{p}\omega_i\bar{x}_i}{\sum_{i=1}^{p}\omega_i\mu_i}; \quad \sum_{i=1}^{p}\omega_i = 1, \text{ [Sahai \textit{et al} (1980)]}$$

$$\hat{\mu}_g^{(5)} = \bar{y}\prod_{i=1}^{p}\left(\frac{\mu_i}{\bar{x}_i}\right)^{\omega_i}, \quad \sum_{i=1}^{p}\omega_i = 1, \text{ [Mohanty and Pattanaik (1984)]}$$

$$\hat{\mu}_g^{(6)} = \bar{y}\left(\sum_{i=1}^{p}\frac{\omega_i\bar{x}_i}{\mu_i}\right)^{-1}, \quad \sum_{i=1}^{p}\omega_i = 1, \text{ [Mohanty and Pattanaik (1984)]}$$

$$\hat{\mu}_g^{(7)} = \bar{y}\prod_{i=1}^{p}\left(\frac{\bar{x}_i}{\mu_i}\right)^{\omega_i}, \quad \sum_{i=1}^{p}\omega_i = 1, \text{ [Tuteja and Bahl (1991)]}$$

$$\hat{\mu}_g^{(8)} = \bar{y}\left[\sum_{i=1}^{p}\frac{\omega_i\mu_i}{\bar{x}_i}\right]^{-1}, \quad \sum_{i=1}^{p}\omega_i = 1, \text{ [Tuteja and Bahl (1991)]}$$

$$\hat{\mu}_g^{(9)} = \bar{y}\left[\omega_{p+1} + \sum_{i=1}^{p}\omega_i\left(\frac{\mu_i}{\bar{x}_i}\right)\right], \quad \sum_{i=1}^{p+1}\omega_i = 1.$$

$$\hat{\mu}_g^{(10)} = \bar{y}\left[\omega_{p+1} + \sum_{i=1}^{p}\omega_i\left(\frac{\bar{x}_i}{\mu_i}\right)\right], \quad \sum_{i=1}^{p+1}\omega_i = 1.$$

$$\hat{\mu}_g^{(11)} = \bar{y}\left[\sum_{i=1}^{q}\omega_i\left(\frac{\mu_i}{x_i}\right) + \sum_{i=q+1}^{p}\left(\frac{\hat{x}_i}{\mu_i}\right)\right]; \quad \left(\sum_{i=1}^{q}\omega_i + \sum_{i=q+1}^{p}\omega_i\right)^{=1}; \text{ [Srivastava (1965) and Rao and Mudhalkar (1967)]}$$

$$\hat{\mu}_g^{(12)} = \bar{y}\prod_{i=1}^{p}\left(\frac{\bar{x}_i}{\mu_i}\right)^{\alpha_i} (\alpha_i\text{'s are suitably constants}) \text{ [Srivastava (1967)]}$$

$$\hat{\mu}_g^{(13)} = \bar{y}\prod_{i=1}^{p}\left\{2 - \left(\frac{\bar{x}_i}{\mu_i}\right)^{\alpha_i}\right\} \text{ [Sahai and Rey (1980)]}$$

$$\hat{\mu}_g^{(14)} = \bar{y}\prod_{i=1}^{p}\frac{\bar{x}_i}{\{\mu_i + \alpha_i(\bar{x}_i - \mu_i)\}} \text{ [Walsh (1970)]}$$

$$\hat{\mu}_g^{(15)} = \bar{y}\exp\left\{\sum_{i=1}^{p}\theta_i \log u_i\right\} \text{ [Srivastava (1971)]}$$

$$\hat{\mu}_g^{(16)} = \bar{y}\exp\left\{\sum_{i=1}^{p}\theta_i(u_i - 1)\right\} \text{ [Srivastava (1971)]}$$

$$\hat{\mu}_g^{(17)} = \bar{y}\sum_{i=1}^{p}\omega_i \exp\{(\theta_i/\omega_i)\log u_i\}, \quad \sum_{i=1}^{p}\omega_i = 1, \text{ [Srivastava (1971)]}$$

$$\hat{\mu}_g^{(18)} = \bar{y} + \sum_{i=1}^{p}\alpha_i(\bar{x}_i - \mu_i)$$

etc. may be identified as particular members of the suggested family of estimators $\hat{\mu}_g$. The MSE of these estimators can be obtained from (2.4).

It is well known that

$$\mathrm{V}(\bar{y}) = (\mu_0^2/n)(C_0^2 + C_{(0)}^2) \tag{2.8}$$

It follows from (2.6) and (2.8) that the minimum variance of $\hat{\mu}_g$ is no longer than conventional unbiased estimator $\bar{y}$.

On substituting $\sigma_{(0)}^2=0$, $\sigma_{(i)}^2=0$ $\forall i=1,2,\ldots,p$ in the equation (2.4), we obtain the no-measurement error case. In that case, the MSE of $\hat{\mu}_g$, is given by

$$\begin{aligned}\mathrm{MSE}(\hat{\mu}_g) &= \frac{1}{n}\left[C_0^2\mu_0^2 + 2\mu_0 b^T g^{*(1)}(\mu_0, e^T) + \left(g^{*(1)}(\mu_0, e^T)\right)^T A^*\left(g^{*(1)}(\mu_0, e^T)\right)\right] \\ &= \mathrm{MSE}(\hat{\mu}_g{}^*)\end{aligned} \tag{2.9}$$

where

$$\begin{aligned}\hat{\mu}_g &= g^*\left(\bar{Y}, \frac{\bar{X}_1}{\mu_1}, \frac{\bar{X}_2}{\mu_2}, \cdots, \frac{\bar{X}_p}{\mu_p}\right) \\ &= g^*(\bar{Y}, U^T)\end{aligned} \tag{2.10}$$

and $\overline{\mathrm{Y}}$ and $\overline{\mathrm{X}}_{\mathrm{i}}\,(i=1,2,\cdots,p)$ are the sample means of the characteristics Y and $X_i$ based on true measurements. ($Y_j$,$X_{ij}$, i=1,2,…,p; j=1,2,…,n). The family of estimators $\hat{\mu}_g*$ at (2.10) is a generalized version of Srivastava (1971, 80).

The MSE of $\hat{\mu}_g*$ is minimized for

$$g*^{(1)}(\mu_0,e^T) = -A*^{-1}\, b\mu_0 \tag{2.11}$$

Thus the resulting minimum MSE of $\hat{\mu}_g*$ is given by

$$\begin{aligned}\text{min.MSE}(\hat{\mu}_g *) &= \frac{{\mu_0}^2}{n}\left[C_0^2 - b^T A*^{-1}\, b\right]\\ &= \frac{{\sigma_0}^2}{n}\left(1-R^2\right)\end{aligned} \tag{2.12}$$

where A*=[$a^*_{ij}$] be a p×p matrix with $a^*_{ij} = \rho_{ij}C_iC_j$ and R stands for the multiple correlation coefficient of Y on $X_1$,$X_2$,…,$X_p$.

From (2.6) and (2.12) the increase in minimum MSE $(\hat{\mu}_g)$ due to measurement errors is obtained as

$$\begin{aligned}\text{min.MSE}(\hat{\mu}_g) - \text{min.MSE}(\hat{\mu}_g *) &= \left(\frac{\mu_0^2}{n}\right)\left[C_{(0)}^2 + b^T A*^{-1}\, b - b^T A^{-1} b\right]\\ &> 0\end{aligned}$$

This is due to the fact that the measurement errors introduce the variances fallible measurements of study variate Y and auxiliary variates $X_i$. Hence there is a need to take the contribution of measurement errors into account.

## 3. BIASES AND MEAN SQUARE ERRORS OF SOME PARTICULAR ESTIMATORS IN THE PRESENCE OF MEASUREMENT ERRORS.

To obtain the bias of the estimator $\hat{\mu}_g$, we further assume that the third partial derivatives of $g(\overline{y},u^T)$ also exist and are continuous and bounded. Then expanding $g(\overline{y},u^T)$ about the point $(\overline{y},u^T) = (\mu_0,e^T)$ in a third-order Taylor's series we obtain

$$\hat{\mu}_g = g(\mu_0,e^T) + (\overline{y}-\mu_0)\left.\frac{\partial g(\cdot)}{\partial \overline{y}}\right|_{(\mu_0,e^T)} + (u-e)^T g^{(1)}(\mu_0,e^T)$$

$$+\frac{1}{2}\left\{(\bar{y}-\mu_0)^2\frac{\partial^2 g(\cdot)}{\partial\bar{y}^2}\bigg|_{(\mu_0,u^T)}+2(\bar{y}-\mu_0)(u-e)^T g^{(1)}(\mu_0,e^T)\right.$$

$$\left.+(u-e)^T\left(g^{(2)}(\mu_0,e^T)\right)(u-e)\right\}$$

$$+\frac{1}{6}\left\{(\bar{y}-\mu_0)\frac{\partial}{\partial\bar{y}}+(u-e)\frac{\partial}{\partial u}\right\}^3 g\left(\bar{y}^*,u^{*T}\right)$$

(3.1)

where $g^{(12)}(\mu_0,e^T)$ denotes the matrix of second partial derivatives of $g\left(\bar{y},u^T\right)$ at the point $\left(\bar{y},u^T\right)=\left(\mu_0,e^T\right)$.

Noting that

$$g\left(u_0 e^T\right)=\mu_0$$

$$\frac{\partial g(\cdot)}{\partial\bar{y}}\bigg|_{(\mu_0,e^T)}=1$$

$$\frac{\partial^2 g(\cdot)}{\partial\bar{y}^2}\bigg|_{(\mu_0,e^T)}=0$$

and taking expectation we obtain the bias of the family of estimators $\hat{\mu}_g$ to the first degree of approximation,

$$B\left(\hat{\mu}_g\right)=\frac{1}{2}\left[E\left\{(u-e)^T\left(g^{(2)}(\mu_0,e^T)\right)(u-e)\right\}+2\left(\frac{\mu_0}{n}\right)b^T g^{(12)}(\mu_0,e^T)\right]$$

(3.2)

where $b^T=(b_1,b_2,\ldots,b_p)$ with $bi=\rho_{oi}C_0C_i$; $(i=1,2,\ldots, p)$. Thus we see that the bias of $\hat{\mu}_g$ depends also upon the second order partial derivatives of the function on $g\left(\bar{y},u^T\right)$ at the point $(\mu_0,e^T)$, and hence will be different for different optimum estimators of the family.

The biases and mean square errors of the estimators $\hat{\mu}_g^{(i)}; i=1\text{ to }18$ up to terms of order $n^{-1}$ along with the values of $g^{(1)}(\mu_0,e^T)$, $g^{(2)}(\mu_0,e^T)$ and $g^{(12)}(\mu_0,e^T)$ are given in the Table 3.1.

**Table 3.1 Biases and mean squared errors of various estimators of $\mu_0$**

| ESTIMATOR | $g^{(1)}(\mu_0,e^T)$ | $g^{(2)}(\mu_0,e^T)$ | $g^{(12)}(\mu_0,e^T)$ | BIAS | MSE |
|---|---|---|---|---|---|
| $\hat{\mu}_g^{(1)}$ | $-\mu_0\underset{\sim}{\omega}$ | $2\mu_0\underset{\sim}{W}_{p\times p}$ where $W_{pxp}=dig(\omega_1,\omega_2,\ldots,\omega_p)$ | $-\underset{\sim}{\omega}$ | $\left(\frac{\mu_0}{n}\right)\left(C^T\underset{\sim}{W}_{p\times p}-b^T\underset{\sim}{\omega}\right)$ where $C^T=\left(C_1^2+C_{(1)}^2,C_2^2+C_{(2)}^2,\cdots,C_p^2+C_{(p)}^2\right)$ | $\left(\frac{\mu_0{}^2}{n}\right)\left[C_0^2+C_{(o)}^2-2b^T\underset{\sim}{\omega}+\underset{\sim}{\omega}^T A\underset{\sim}{\omega}\right]$ |
| $\hat{\mu}_g^{(2)}$ | $\mu_0\underset{\sim}{\omega}$ | $\underset{\sim}{O}_{p\times p}$ (null matrix) | $\underset{\sim}{\omega}$ | $\left(\frac{\mu_0}{n}\right)b^T\underset{\sim}{\omega}$ | $\left(\frac{\mu_0^2}{n}\right)\left[C_0^2+C_{(0)}^2+2b^T\underset{\sim}{\omega}+\underset{\sim}{\omega}^T A\underset{\sim}{\omega}\right]$ |
| $\hat{\mu}_g^{(3)}$ | $-\frac{\mu_0\underset{\sim}{\omega}^*}{\underset{\sim}{\omega}^T\underset{\sim}{\mu}}$ where $\underset{\sim}{\omega}^{*T}=$ $(\omega_1{}^*,\omega_2{}^*,\ldots,\omega_p{}^*)$ with $(\omega_i{}^*=\omega_i\,\mu_i)$ $(i=1,2,\ldots,p)$ | $\frac{2\mu_0\underset{\sim}{\omega}^*\underset{\sim}{\omega}^{*T}}{\underset{\sim}{\omega}^T\underset{\sim}{\mu}\,\underset{\sim}{\mu}^T\underset{\sim}{\omega}}$ | $-\frac{\underset{\sim}{\omega}^*}{\underset{\sim}{\omega}^T\underset{\sim}{\mu}}$ | $\left(\frac{\mu_0}{n}\right)\left(\frac{\underset{\sim}{\omega}^{*T}A\underset{\sim}{\omega}}{\underset{\sim}{\omega}^T\underset{\sim}{\mu}\,\underset{\sim}{\mu}^T\underset{\sim}{\omega}}-\frac{b^T\underset{\sim}{\omega}^*}{\underset{\sim}{\omega}^T\underset{\sim}{\mu}}\right)$ | $\left(\frac{\mu_0^2}{n}\right)\left[C_0^2+C_{(0)}^2-\frac{2b^T\underset{\sim}{\omega}^*}{\underset{\sim}{\omega}^T\underset{\sim}{\mu}}+\frac{\underset{\sim}{\omega}^{*T}A\underset{\sim}{\omega}^*}{\underset{\sim}{\omega}^T\underset{\sim}{\mu}\,\underset{\sim}{\mu}^T\underset{\sim}{\omega}}\right]$ |
| $\hat{\mu}_g^{(4)}$ | $\frac{\mu_0\underset{\sim}{\omega}}{\underset{\sim}{\omega}^T\underset{\sim}{\mu}}$ | $\underset{\sim}{O}_{p\times p}$ (null matrix) | $\frac{\underset{\sim}{\omega}}{\underset{\sim}{\omega}^T\underset{\sim}{\mu}}$ | $\left(\frac{\mu_0}{n}\right)\frac{b^T\underset{\sim}{\omega}}{\underset{\sim}{\omega}^T\underset{\sim}{\mu}}$ | $\left(\frac{\mu_0^2}{n}\right)\left[C_0^2+C_{(0)}^2+\frac{2b^T\underset{\sim}{\omega}}{\underset{\sim}{\omega}^T\underset{\sim}{\mu}}+\frac{\underset{\sim}{\omega}^T A\underset{\sim}{\omega}}{\underset{\sim}{\omega}^T\underset{\sim}{\mu}\,\underset{\sim}{\mu}^T\underset{\sim}{\omega}}\right]$ |
| $\hat{\mu}_g^{(5)}$ | $-\mu_0\underset{\sim}{\omega}$ | $\mu_0\left(\underset{\sim}{\omega}\,\underset{\sim}{\omega}^T+\underset{\sim}{W}_{p\times p}\right)$ | $-\underset{\sim}{\omega}$ | $\left(\frac{\mu_0}{2n}\right)\left[\underset{\sim}{\omega}^T A\underset{\sim}{\omega}+C^T\underset{\sim}{W}_{p\times p}-2b^T\underset{\sim}{\omega}\right]$ | $\left(\frac{\mu_0{}^2}{n}\right)\left[C_0^2+C_{(o)}^2-2b^T\underset{\sim}{\omega}+\underset{\sim}{\omega}^T A\underset{\sim}{\omega}\right]$ |
| $\hat{\mu}_g^{(6)}$ | $-\mu_0 \underset{\sim}{\omega}$ | $2\mu_0 \underset{\sim}{\omega}^T \underset{\sim}{\omega}$ | $-\underset{\sim}{\omega}$ | $\left(\frac{\mu_0}{n}\right)\left[\underset{\sim}{\omega}^T A \underset{\sim}{\omega} - b^T \underset{\sim}{\omega}\right]$ | $\left(\frac{\mu_0^2}{n}\right)\left[C_0^2 + C_{(o)}^2 - 2b^T \underset{\sim}{\omega} + \underset{\sim}{\omega}^T A \underset{\sim}{\omega}\right]$ |
| $\hat{\mu}_g^{(7)}$ | $\mu_0 \underset{\sim}{\omega}$ | $\mu_0\left(\underset{\sim}{\omega}\underset{\sim}{\omega}^T - \underset{\sim}{W}_{p\times p}\right)$ | $\underset{\sim}{\omega}$ | $\left(\frac{\mu_0}{2n}\right)\left[\underset{\sim}{\omega}^T A \underset{\sim}{\omega} - C^T \underset{\sim}{W}_{p\times p} + 2b^T \underset{\sim}{\omega}\right]$ | $\left(\frac{\mu_0^2}{n}\right)\left[C_0^2 + C_{(o)}^2 + 2b^T \underset{\sim}{\omega} + \underset{\sim}{\omega}^T A \underset{\sim}{\omega}\right]$ |
| $\hat{\mu}_g^{(8)}$ | $\mu_0 \underset{\sim}{\omega}$ | $2\mu_0\left(\underset{\sim}{\omega}\underset{\sim}{\omega}^T - \underset{\sim}{W}_{p\times p}\right)$ | $\underset{\sim}{\omega}$ | $\left(\frac{\mu_0}{n}\right)\left[\underset{\sim}{\omega}^T A \underset{\sim}{\omega} - C^T \underset{\sim}{W}_{p\times p} + b^T \underset{\sim}{\omega}\right]$ | $\left(\frac{\mu_0^2}{n}\right)\left[C_0^2 + C_{(o)}^2 + 2b^T \underset{\sim}{\omega} + \underset{\sim}{\omega}^T A \underset{\sim}{\omega}\right]$ |

**Table 3.1 Biases and mean squared errors of various estimators of $\mu_0$**

<table>
<tr><th>ESTIMATOR</th><th>$g^{(1)}(\mu_0,e^T)$</th><th>$g^{(2)}(\mu_0,e^T)$</th><th>$g^{(12)}(\mu_0,e^T)$</th><th>BIAS</th><th>MSE</th></tr>
<tr><td>$\hat{\mu}_g^{(9)}$</td><td>$-\mu_0 \underset{\sim}{\omega}$</td><td>$2\mu_0 \underset{\sim}{W}_{p\times p}$</td><td>$-\underset{\sim}{\omega}$</td><td>$\left(\frac{\mu_0}{n}\right)\left(C^T \underset{\sim}{W}_{p\times p} - b^T \underset{\sim}{\omega}\right)$</td><td>$\left(\frac{\mu_0^2}{n}\right)\left[C_0^2 + C_{(0)}^2 + 2b^T \underset{\sim}{\omega} + \underset{\sim}{\omega}^T A \underset{\sim}{\omega}\right]$</td></tr>
<tr><td>$\hat{\mu}_g^{(10)}$</td><td>$\mu_0 \underset{\sim}{\omega}$</td><td>$\underset{\sim}{O}$</td><td>$\underset{\sim}{\omega}$</td><td>$\left(\frac{\mu_0}{n}\right) b^T \underset{\sim}{\omega}$</td><td>$\left(\frac{\mu_0^2}{n}\right)\left[C_0^2 + C_{(0)}^2 + 2b^T \underset{\sim}{\omega} + \underset{\sim}{\omega}^T A \underset{\sim}{\omega}\right]$</td></tr>
<tr><td rowspan="2">$\hat{\mu}_g^{(11)}$</td><td>$\underset{\sim}{\omega}_{(1)} \mu_0$</td><td>$2\underset{\sim}{W}_{(1)p\times p} \mu_0$</td><td>$\underset{\sim}{\omega}_{(1)}$</td><td rowspan="2">$\left(\frac{\mu_0}{n}\right)\left(C^{*T} \underset{\sim}{W}_{(1)} - b^T \underset{\sim}{\omega}_{(1)}\right)$<br><br>$C^{*T}=(C_1^2+C_{(1)}^2,\ldots,$<br>$C_q^2+C_{(q)}^2;\ldots 0)$</td><td rowspan="2">$\left(\frac{\mu_0^2}{n}\right)\left[C_0^2 + C_{(o)}^2 - 2b^T \underset{\sim}{\omega}_{(1)} + \underset{\sim}{\omega}^T_{(1)} A \underset{\sim}{\omega}_{(1)}\right]$</td></tr>
<tr><td colspan="3">where $\underset{\sim}{\omega}_{(1)}=(-\omega_1,-\omega_2,\ldots,-\omega_q,-\omega_{q+1},\ldots,\omega_p)_{1\times p}$<br><br>$\underset{\sim}{W}_{(1)p\times p}=\begin{bmatrix}\omega_1 & 0 & 0 & \cdots & 0 & 0 & \cdots & 0\\ 0 & \omega_2 & 0 & \cdots & 0 & 0 & \cdots & 0\\ 0 & 0 & \omega_3 & \cdots & 0 & 0 & \cdots & 0\\ \vdots & \vdots & \vdots & \cdots & \vdots & \vdots & \cdots & \vdots\\ 0 & 0 & 0 & \cdots & \omega_q & 0 & \cdots & 0\\ 0 & 0 & 0 & \cdots & 0 & 0 & \cdots & 0\\ \vdots & \vdots & \vdots & \cdots & \vdots & \vdots & \cdots & \vdots\\ 0 & 0 & 0 & \cdots & 0 & 0 & \cdots & 0\end{bmatrix}_{p\times p},$</td></tr>
</table>

| | | | | | |
|---|---|---|---|---|---|
| $\hat{\mu}_g^{(12)}$ | $\underset{\sim}{\alpha}\,\mu_0$ | $\mu_0\left(\underset{\sim}{\alpha}\,\underset{\sim}{\alpha}^T - \underset{\sim}{\propto}_{p\times p}\right)$ | $-\underset{\sim}{\alpha}$ | $\left(\frac{\mu_0}{2n}\right)\left[\underset{\sim}{\alpha}^T A\underset{\sim}{\alpha} - C^T\underset{\sim}{\propto}_{p\times p} + 2b^T\underset{\sim}{\alpha}\right]$ | $\left(\frac{\mu_0{}^2}{n}\right)\left[C_0^2 + C_{(o)}^2 + 2b^T\underset{\sim}{\alpha} + \underset{\sim}{\alpha}^T A\,\underset{\sim}{\alpha}\right]$ |
| | where $\underset{\sim}{\alpha}^T = (\alpha_1,\alpha_2,\ldots,\alpha_p)_{1\times p}$ $\underset{\sim}{\propto} = \mathrm{diag}(\alpha_1,\alpha_2,\ldots,\alpha_p)$ | | | | |
| $\hat{\mu}_g^{(13)}$ | $-\underset{\sim}{\alpha}\,\mu_0$ | $-\mu_0\left(\underset{\sim}{\alpha}\,\underset{\sim}{\alpha}^T - \underset{\sim}{\propto}_{p\times p}\right)$ | $-\underset{\sim}{\alpha}$ | $\left(\frac{\mu_0}{2n}\right)\left[C^T\underset{\sim}{\propto}_{p\times p} - \underset{\sim}{\alpha}^T A\underset{\sim}{\alpha} - 2b^T\underset{\sim}{\alpha}\right]$ | $\left(\frac{\mu_0{}^2}{n}\right)\left[C_0^2 + C_{(o)}^2 - 2b^T\underset{\sim}{\alpha} + \underset{\sim}{\alpha}^T A\,\underset{\sim}{\alpha}\right]$ |
| $\hat{\mu}_g^{(14)}$ | $-\underset{\sim}{\alpha}\,\mu_0$ | $2\mu_0\,\underset{\sim}{\alpha}\,\underset{\sim}{\alpha}$ | $-\underset{\sim}{\alpha}$ | $\left(\frac{\mu_0}{n}\right)\left[\underset{\sim}{\alpha}^T A\underset{\sim}{\alpha} - C^T - b^T\underset{\sim}{\alpha}\right]$ | $\left(\frac{\mu_0{}^2}{n}\right)\left[C_0^2 + C_{(o)}^2 - 2b^T\underset{\sim}{\alpha} + \underset{\sim}{\alpha}^T A\,\underset{\sim}{\alpha}\right]$ |

**Table 3.1 Biases and mean squared errors of various estimators of $\mu_0$**

| ESTIMATOR | $g^{(1)}(\mu_0,e^T)$ | $g^{(2)}(\mu_0,e^T)$ | $g^{(12)}(\mu_0,e^T)$ | BIAS | MSE |
|---|---|---|---|---|---|
| $\hat{\mu}_g^{(15)}$ | $\mu_0\underset{\sim}{\theta}$,<br>where $\underset{\sim}{\Theta}_{p\times p}$ =diag{$\theta_1,\theta_2,\ldots\ \theta_p$} | $\mu_0\left(\underset{\sim}{\theta}\underset{\sim}{\theta}^T-\underset{\sim}{\Theta}_{p\times p}\right)$, | $v\,0\ \underset{\sim}{\theta}$ | $\left(\frac{\mu_0}{2n}\right)\left[\underset{\sim}{\theta}^T A\underset{\sim}{\theta}-C^T\underset{\sim}{\Theta}_{p\times p}+2b^T\underset{\sim}{\theta}\right]$ | $\left(\frac{\mu_0^2}{n}\right)\left[C_0^2+C_{(0)}^2+2b^T\underset{\sim}{\theta}+\underset{\sim}{\theta}^T A\underset{\sim}{\theta}\right]$ |
| $\hat{\mu}_g^{(16)}$ | $\mu_0\underset{\sim}{\theta}$ | $\mu_0\underset{\sim}{\theta}\underset{\sim}{\theta}^T$ | $\underset{\sim}{\theta}$ | $\left(\frac{\mu_0}{2n}\right)\left[\underset{\sim}{\theta}^T A\underset{\sim}{\theta}+2b^T\underset{\sim}{\theta}\right]$ | $\left(\frac{\mu_0^2}{n}\right)\left[C_0^2+C_{(0)}^2+2b^T\underset{\sim}{\theta}+\underset{\sim}{\theta}^T A\underset{\sim}{\theta}\right]$ |
| $\hat{\mu}_g^{(17)}$ | $\mu_0\underset{\sim}{\theta}$<br>where $\underset{\sim}{\Theta}*_{p\times p}$ =diag{$\theta_1\left(\frac{\theta_1}{\omega_1}-1\right)\ldots,\theta_p\left(\frac{\theta_p}{\omega_p}-1\right)$} | $\underset{\sim}{\Theta}*_{p\times p}\ \mu_0$, | $\underset{\sim}{\theta}$ | $\left(\frac{\mu_0}{2n}\right)\left[C^T\underset{\sim}{\Theta}*_{p\times p}+2b^T\underset{\sim}{\theta}\right]$ | $\left(\frac{\mu_0{}^2}{n}\right)\left[C_0^2+C_{(o)}^2+2b^T\underset{\sim}{\theta}+\underset{\sim}{\theta}^T A\underset{\sim}{\theta}\right]$ |
| $\hat{\mu}_g^{(18)}$ | $\underset{\sim}{\alpha}*$<br>where $\underset{\sim}{\alpha}*^T$ =($\alpha_{1*},\alpha_{2*},\ldots,\alpha_{p*}$) with $\underset{\sim}{\alpha}*_i$ =($\alpha_i,\mu_i$, i=1,2,...,p) | $\underset{\sim}{O}_{p\times p}$ | $\underset{\sim}{O}_{p\times p}$ | Unbiased | $\left(\frac{1}{n}\right)\left[C_0^2+C_{(o)}^2+2\mu_0 b^T\underset{\sim}{\alpha}*+\underset{\sim}{\alpha}*^T A\underset{\sim}{\alpha}*\right]$ |

## 4. ESTIMATORS BASED ON ESTIMATED OPTIMUM

It may be noted that the minimum MSE (2.6) is obtained only when the optimum values of constants involved in the estimator, which are functions of the unknown population parameters $\mu_0$, b and A, are known quite accurately.

To use such estimators in practice, one has to use some guessed values of the parameters $\mu_0$, b and A, either through past experience or through a pilot sample survey. Das and Tripathi (1978, sec.3) have illustrated that even if the values of the parameters used in the estimator are not exactly equal to their optimum values as given by (2.5) but are close enough, the resulting estimator will be better than the conventional unbiased estimator $\bar{y}$. For further discussion on this issue, the reader is referred to Murthy (1967), Reddy (1973), Srivenkataramana and Tracy (1984) and Sahai and Sahai (1985).

On the other hand if the experimenter is unable to guess the values of population parameters due to lack of experience, it is advisable to replace the unknown population parameters by their consistent estimators. Let $\hat{\phi}$ be a consistent estimator of $\phi=A^{-1}b$. We then replace $\phi$ by $\hat{\phi}$ and also $\mu_0$ by $\bar{y}$ if necessary, in the optimum $\hat{\mu}_g$ resulting in the estimator $\hat{\mu}_{g(est)}$, say, which will now be a function of $\bar{y}$, u and $\phi$. Thus we define a family of estimators (based on estimated optimum values) of $\mu_0$ as

$$\hat{\mu}_{g(est)} = g^{**}\left(\bar{y}, u^T, \hat{\phi}^T\right) \tag{4.1}$$

where g**(·) is a function of $\left(\bar{y}, u^T, \hat{\phi}^T\right)$ such that

$$g^{**}\left(\mu_0, e^T, \phi^T\right) = \mu_0 \text{ for all } \mu_0,$$

$$\Rightarrow \left.\frac{\partial g^{**}(\cdot)}{\partial \bar{y}}\right|_{\left(\mu_0, e^T, \phi^T\right)} = 1$$

$$\left.\frac{\partial g^{**}(\cdot)}{\partial u}\right|_{\left(\mu_0, e^T \phi^T\right)} = \left.\frac{\partial g(\cdot)}{\partial u}\right|_{(\mu_0, e)^T} = -\mu_0 A^{-1} b = -\mu_0 \phi \tag{4.2}$$

and

$$\left.\frac{\partial g^{**}(\cdot)}{\partial \hat{\phi}}\right|_{\left(\mu_0, e^T, \phi^T\right)} = 0$$

With these conditions and following Srivastava and Jhajj (1983), it can be shown to the first degree of approximation that

$$\mathrm{MSE}\left(\hat{\mu}_{g(est)}\right) = \min.\mathrm{MSE}\left(\hat{\mu}_{g}\right)$$
$$= \left(\frac{\mu_0^{\ 2}}{n}\right)\left[C_0^2 + C_{(0)}^2 - b^T A^{-1} b\right]$$

Thus if the optimum values of constants involved in the estimator are replaced by their consistent estimators and conditions (4.2) hold true, the resulting estimator $\hat{\mu}_{g(est)}$ will have the same asymptotic mean square error, as that of optimum $\hat{\mu}_g$. Our work needs to be extended and future research will explore the computational aspects of the proposed algorithm.